\documentclass[conference]{IEEEtran}
\IEEEoverridecommandlockouts

\usepackage{graphics}
\usepackage{graphicx} 
\usepackage{rotating}
\usepackage{dcolumn}
\usepackage{hyperref}
\usepackage{longtable}
\usepackage{multirow}
\usepackage{amsmath}
\usepackage{amssymb}
\usepackage{xspace}
\usepackage{textcase}
\usepackage{widetext}
\usepackage{float}
\usepackage{amsthm}
\usepackage{multicol}
\usepackage{standalone}
\usepackage{tikz}
\usepackage{tikzscale}
\usepackage{pgfplots}
\pgfplotsset{compat=1.18}
\usetikzlibrary{trees}
\usetikzlibrary{arrows.meta}
\usetikzlibrary{positioning}
\usetikzlibrary{arrows, math, angles, quotes, shapes.geometric, calc}
\usepackage{commath}
\usepackage{xcolor}
\usepackage{hyperref}
\usepackage{subfiles}

\def\BibTeX{{\rm B\kern-.05em{\sc i\kern-.025em b}\kern-.08em
    T\kern-.1667em\lower.7ex\hbox{E}\kern-.125emX}}

\begin{document}

\title{Time Scale Separation and Hierarchical Control with the Koopman Operator\
\thanks{This work was supported by the DOE Wind Energy Technologies Office.}
}

\author{\IEEEauthorblockN{1\textsuperscript{st} Craig Bakker}
\IEEEauthorblockA{\textit{National Security Directorate} \\
\textit{Pacific Northwest National Laboratory}\\
Richland, Washington \\
craig.bakker@pnnl.gov}
}

\maketitle


\begin{abstract}
The Koopman Operator (KO) is a mathematical construct that maps nonlinear (state space) dynamics to corresponding linear dynamics in an infinite-dimensional functional space.  For practical applications, finite-dimensional approximations can be constructed with machine learning.  The linearity of the KO facilitates the use of linear tools and theories on nonlinear dynamical systems without relying on local approximations.  Additionally, the KO has the ability to incorporate forms of domain knowledge such as stability and known variable interactions.  It therefore constitutes part of the broader interest in domain-aware or physics-informed machine learning.  Hierarchical control and time scale separation are two common properties in engineered systems -- often found together -- that can pose both computational and practical challenges.  These properties provide opportunities for further KO-based domain knowledge incorporation, but this has seldom been taken advantage of in the KO literature.  This paper focuses on developing and using domain-aware Koopman formulations for systems with hierarchical control, systems with time scale separation, and systems with both.  We show how those formulations can be leveraged to quantify the impact of cross-scale interactions on overall stability and to calculate optimal control policies at both fast and slow time scales.
\end{abstract}


\section{Introduction}

The Koopman Operator (KO) is a mathematical construct that maps nonlinear (state space) dynamics to corresponding linear dynamics in an infinite-dimensional functional space of observables \cite{budisic12jsr}.  In practice, the typical approach is to compute a discrete-time, finite-dimensional approximation of the KO, and there are various ways of doing this \cite{bakker2020koopman}.  For example, if we have dynamics $\dot{x} = f \left(x\right)$, our (discrete-time) KO model would be $\psi_{x,t+1} = K \psi_{x,t}$, where $\psi_{x,t} = \psi_x \left(x_t\right)$; $\psi_{x,t}$ is a vector of Koopman observables, and the function $\psi_x \left(x\right)$ is said to ``lift'' the states $x$ from the state space to the space of observables.  Once an approximation has been calculated, it can then be used for control \cite{huang2018feedback,korda2018linear,king2021solving} and analysis purposes \cite{sinha2022data,mezic15cp,surana2016linear}.  The linearity of the KO facilitates the use of linear tools and theories on nonlinear dynamical systems without relying on local approximations.  This can be particularly useful for optimal control -- both Model Predictive Control (MPC) and Koopman forms of the Linear Quadratic Regulator (LQR).  It can also be useful (or even necessary) for use in adversarial or game-theoretic contexts: computing optimal sensor bias attacks \cite{bakker2022deception}, solving zero-sum differential games \cite{bakker2025operator}, or solving defender-attacker-defender games \cite{oster2023multi}.

The KO also has the ability to incorporate forms of domain knowledge such as stability \cite{king2021solving}, symmetry \cite{sinha2020koopman}, and known variable interactions \cite{bakker2022deception}; it therefore constitutes part of the broader interest in domain-aware or physics-informed machine learning \cite{alber2019integrating}.  Hierarchical control and time scale separation are two common properties in engineered systems -- often found together -- that can can pose both computational and practical challenges.  They also provide opportunities for further domain knowledge incorporation, though.  Bakker et al. \cite{bakker2022deception} present a hierarchical control formulation of the KO that recognizes the structure of the interactions between supervisory control inputs, actuators, and system states, but the formulation is not really generalized beyond the application at hand in that paper.  Similarly, Champion et al. \cite{champion2019discovery} consider Koopman learning for systems with multiple time scales, but the multi-scale \emph{discrete-time}  Koopman formulation is not explicitly presented, and thus it does not appear that the authors leverage the separation of scales in their construction of the Koopman model (though they do discuss the use of time scale separation in the \emph{training} of the Koopman model).

This paper focuses on developing and using domain-aware Koopman formulations for systems with hierarchical control, systems with time scale separation, and systems with both.  The combination of the two is motivated by the fact that the lower-level controls in a hierarchical system often operate at a faster time scale than the supervisory controls.  We can then leverage these formulations to quantify the impact of cross-scale interactions on stability in various ways.  The paper also shows how to replace the existing lower-level Proportional-Integral (PI) control with a Koopman-based LQR control process, both with and without time scale separation, and compares their respective performances.

\section{Koopman Model Description}

\subsection{Time Scale Separation}

Let us consider an uncontrolled time scale separated system that, per the standard singular perturbation model of Khalil \cite[Ch. 11]{khalil2002nonlinear}, converges to a Differential Algebraic Equation (DAE) as the scales separate (i.e., as $\epsilon \rightarrow 0$):

\begin{gather}
		\begin{array}{c}\dot{x}  = f \left(x,y\right) \\
		\dot{y}  = \frac{1}{\epsilon} g \left(x,y\right) \end{array} \xrightarrow{\epsilon \rightarrow 0} \begin{array}{c} \dot{x}  = f \left(x,y\right) \\
		0  = g \left(x,y\right) \end{array}
\end{gather}

We can then produce a discrete-time Koopman representation that accounts for the separation of time scales:

\begin{align}
		\psi_x \left(x_{t+\Delta t}\right) &= K_{xx} \left(\psi_x \left(x_t\right) - K_{xy} \frac{1}{m} \sum \limits_{n=0}^{m-1} \psi_y \left(y_{t+n\tau}\right) \right) \nonumber \\
		&+ K_{xy} \frac{1}{m} \sum \limits_{n=0}^{m-1} \psi_y \left(y_{t+n\tau}\right) \label{eq:tss slow} \\
		\psi_y \left(y_{t+\left(n+1\right)\tau}\right) &= K_{yy} \left(\psi_y \left(y_{t+n\tau}\right) - K_{yx} \psi_x \left(x_t\right)\right) \nonumber \\
		&+ K_{yx} \psi_x \left(x_t\right), n \in \left\{0,1,\ldots, m-1\right\} \label{eq:tss fast}
\end{align}

\noindent where $\Delta t = m \tau$, $\tau \propto \epsilon$ captures the separation of scales and $\psi_{\left(\cdot\right)}$ are our observables (split by time scale).  To capture the time scale separation, the fast dynamics in (\ref{eq:tss fast}) evolve $m$ steps for every step of the slow dynamics while the slow variables $x$ are essentially treated as constant inputs.  Here and in the rest of the paper, we use and modify the stability-assuming (or enforcing) KO formulation of King et al. \cite{king2021solving}.  Both this and the averaging over $\psi_y \left(t_{t + n \tau}\right)$ in (\ref{eq:tss slow}) facilitate later calculations.  If the fast dynamics are stable, the $\epsilon \rightarrow 0$ limit produces

\begin{gather}
		\lim_{n \rightarrow \infty} \psi_y \left(y_{t+n\tau}\right) = \psi_y \left(y^* \left(x_t\right)\right) = K_{yx} \psi_x \left(x_t\right) \label{eq:y lim 1} \\
		\lim_{m \rightarrow \infty}  K_{xy} \frac{1}{m} \sum \limits_{n=0}^{m-1} \psi_y \left(y_{t+n\tau}\right) = K_{xy} \psi_y \left(y^* \left(x_t\right)\right) \label{eq:y lim 2}
\end{gather}

\noindent where $0 = g \left(x_t,y^* \left(x_t\right)\right)$ \cite{king2021solving,khalil2002nonlinear}; as the scales separate, the behavior of $y$ becomes more fully characterized by the equilibrium $y^* \left(x_t\right)$, which we can calculate, in terms of the observables, from (\ref{eq:tss fast}).  The slow dynamics then become

\begin{align}
		\psi_x \left(x_{t + \Delta t} \right) &= K_{comb,xx} \psi_x \left(x_t\right) \\
		K_{comb,xx} &= K_{xx} + \left(I - K_{xx}\right) K_{xy} K_{yx}
\end{align}

\subsection{Hierarchical Control}
\label{sec:hierarchical control}

Let us define a hierarchically controlled system as 

\begin{align}
		\dot{y} &= f \left(y,w\right) \\
		\dot{w} &= h \left(y,w,u\right) \\
		J \left(u\right) &= \int \limits_0^T L \left(y,w,u\right) dt
\end{align}

\noindent where $u$ is a set of control inputs, $w$ are actuator variables, $y$ are states, and $J$ is the running cost.  Adding a terminal cost is straightforward but lengthens the derivations, so we omit it here.  Note that the control inputs can only affect the states indirectly, via the actuators, and in anticipation of the next section, we treat the control inputs as constant (relative to $y$ and $w$).  We can produce a discrete-time Koopman representation that captures the variable dependencies:

\begin{align}
		\psi_y \left(y_{t+1}\right) &= K_{yy} \left(\psi_y \left(y_t\right) - K_{yw} \psi_w \left(w_t\right) \right) \nonumber \\
		&+ K_{yw} \psi_w \left(w_t\right) 
		\label{eq:hierarchical dyn 1} \\
		\psi_w \left(w_{t+1}\right) &= K_{ww} \left(\psi_w \left(w_t\right) - K_{wy} \psi_y \left(y_t\right) \right. \nonumber \\
		& \left. - K_{wu} \psi_u \left(u\right) \right) + K_{wy} \psi_y \left(y_t\right) + K_{wu} \psi_u \left(u\right)
		\label{eq:hierarchical dyn 2}
\end{align}

There are additional ways to incorporate domain knowledge into the structure of the $K$ matrices \cite{bakker2022deception}.  We also define a quadratic (in $\psi_{\left(\cdot\right)}$) discrete-time cost approximation cost $\hat{L}$

\begin{align}
		\hat{L} &= \sum_t \left[\sum_{i,j} \psi_{i,t}^T Q_{ij} \psi_{j,t} + \sum_i c_i^T \psi_{i,t} + c_0 \right]
		\label{eq:cost approx}
\end{align}

\noindent where $i,j \in \left\{u,w,y\right\}$ and $\psi_{u,t} = \psi_u$ is constant in time; we assume that $Q_{ii}$ is symmetric.  We could learn $\hat{L}$ from data, but if $L$ is known, it may be easier to specify $\hat{L}$ as a Riemann approximation to $L$; e.g., if $L = y^2$ and $\psi_y$ is state-inclusive, we would specify the first diagonal entry of $Q_{yy}$ to be $\Delta t$ with all other $Q$ and $c$ components being zero.  Alternatively, we control the actuators $w$ with an optimal feedback policy.  Let us assume an infinite-horizon LQR-like policy, for a fixed control input $u$, where our value function has both linear and quadratic terms:

\begin{align}
		V \left(y\right) &= \psi_y^T P \psi_y + p^T \psi_y \\
		\psi_{w,t} &= -F \psi_{y,t} - d
		\label{eq:lqr feedback 1}
\end{align}

Both $p$ and $d$ will be functions of $u$.  We can then use Bellman's equation to solve for $P$, $F$, $d$, and $p$, thereby producing our feedback control policy.  The full equations for the case with time scale separation and hierarchical control are presented in the following section; the solution without time scale separation simply involves omitting any terms containing a subscript $x$.  Also, if LQR control is being substituted for PI control, the LQR objective may not be identical to the overall cost -- e.g, the LQR control may have a cost term of the form $\left( y - u\right)^2$ to mimic the set point matching functionality of PI control -- but let us assume that we can still represent it in the form of (\ref{eq:cost approx}).  This is only possible because of the structure we have imposed upon the Koopman representation of the system: if we did not distinguish between $w$ and $y$, for example, and instead combined them in a single set of observable functions, we would not be able to make this swap.  

A problem we may run into is that $\psi_w$ is not a free vector -- not all values in $\mathbb{R}^{n_{\psi,w}}$ can be achieved by freely varying $w$.  One option would be to extract the optimal $w$ from the calculated $\psi_w$ policy (e.g., from a state-inclusive $\psi_w$) and use that.  Another option would be to enforce the necessary consistency in the training process (by incorporating the LQR calculations into that process).  This would increase training complexity but also potentially improve the quality of the LQR solution.  A third option would be to use $\psi_w \left(w\right) = w$.  The first option produced large inconsistencies, and in the interest of generality, we did not want to assume that $\psi_w \left(w\right) = w$, so we chose the second option.  We thus have two optimal control problems with the objective in (\ref{eq:cost approx}) and any constraints we may wish to employ: one with the original Koopman dynamics in (\ref{eq:hierarchical dyn 1})-(\ref{eq:hierarchical dyn 2}) and one with dynamics governed by (\ref{eq:hierarchical dyn 1}) and (\ref{eq:lqr feedback 1}).

\subsection{Hierarchical Control with Time Scale Separation}

We can now combine the two formulations above for a system where the slow dynamics include the states $x$ and the control inputs $u$, the fast dynamics include the states $y$ and the actuator positions $w$, and the cost function includes all of the variables:

\begin{gather}
		\begin{array}{c}\dot{x} = f \left(x,y,w\right) \\
		\dot{y} = \frac{1}{\epsilon} g \left(x,y,w\right) \\
		\dot{w} = \frac{1}{\epsilon} h \left(x,y,w,u\right) \end{array} \xrightarrow{\epsilon \rightarrow 0} \begin{array}{c} \dot{x} = f \left(x,y,w\right) \\
		0 = g \left(x,y,w\right) \\
		0 = h \left(x,y,w,u\right) \end{array} \\
		J = \int \limits_0^T L \left(w,y,u,x\right) dt
\end{gather}

The discrete-time Koopman representation is presented on the following page in (\ref{eq:comb x dyn})-(\ref{eq:comb w dyn}).  The running cost approximation for these dynamics is 

\begin{align}
		\hat{L} &= \sum_t \frac{1}{m} \sum_{n = 0}^{m-1} C_{t + n \tau} \\
		C_{t + n \tau} &= \sum_{i,j} \psi_{i,t + n\tau}^T Q_{ij} \psi_{j,t+ n \tau} \nonumber \\
		&+ \sum_i c_i^T \psi_{i,t + n \tau} + c_0 \\
		\lim_{\epsilon \rightarrow 0} \sum_{n=0}^{m-1} \frac{C_{t + n \tau}}{m} &= \sum_{i,j} \psi_{i,t}^T Q_{ij} \psi_{j,t} + \sum_i c_i^T \psi_{i,t} + c_0 \\
		\psi_{i,t + n \tau} & = \psi_{i,t} \equiv \psi_{i} \left(i_t\right) \ i \in \left\{x,u\right\} \\
		\psi_{i,t + n \tau} & \equiv \psi_i \left(i_{t + n \tau} \right) \ i \in \left\{w,y\right\} \\
		\psi_{i,t} & \equiv \psi_i \left(i^* \left(x_t,u_t\right)\right)  \ i \in \left\{w,y\right\}
\end{align}

The LQR solution can be obtained by solving Bellman's equation and is shown on the following page in (\ref{eq:comb bellman})-(\ref{eq:comb lqr end}).  The solution $F = Q^{-1} R$ and $d = Q^{-1} q$ relies on the following definitions:

\begin{align}
		Q &= 2 Q_{ww} + 2 K_{yw}^T \left(I - K_{yy}\right)^T P \left(I - K_{yy}\right) K_{yw} \\
		R &= \left[Q_{wy} + Q_{yw}^T + 2 K_{yw}^T \left(I - K_{yy} \right)^T P K_{yy} \right] \psi_{y,t} \\
		q &= \left(Q_{wu} + Q_{uw}^T \right) \psi_{u,t} + c_w + K_{yw}^T \left(I - K_{yy}\right)^T p \nonumber \\
		&+ \left[2 K_{yw}^T \left(I - K_{yy}\right)^T P \left(I - K_{yy} \right) K_{yx} \right. \nonumber \\
		&\left. + Q_{wx} + Q_{xw}^T\right] \psi_{x,t}
\end{align}

Note that $d$ and $p$ are now potentially functions of both $x$ and $u$.  We can solve for the equilibrium of the fast timescale dynamics to produce slow-scale dynamics of the form

\begin{align}
		\psi_w \left(w^* \left(x_t,u_t\right) \right) &= B_{wx} \psi_x \left(x_t\right) + B_{wu} \psi_u \left(u_t\right) \label{eq:comb dyn 1}\\
		\psi_y \left(y^* \left(x_t,u_t\right) \right) &= B_{yx} \psi_x \left(x_t\right) + B_{yu} \psi_u \left(u_t\right) \\
		\psi_x \left(x_{t+1}\right) &= B_{xx} \psi_x \left(x_t\right) + B_{xu} \psi_u \left(u_t\right) \label{eq:comb dyn 3} \\
		0 &= g \left(x,y^* \left(x,u\right),w^* \left(x,u\right)\right) \\
		0 &= h \left(x,y^* \left(x,u\right),w^* \left(x,u\right),u\right) 
\end{align}

\noindent under the $\epsilon \rightarrow 0$ limit; these simplify the cost approximation $\hat{L}$ to be purely a function of $\psi_{x,t}$ and $\psi_{u,t}$.  We can replace the PI control with LQR control, as in the previous section by solving the LQR only at the fast time scale (i.e., assuming $x$ and $u$ are essentially constant) and substituting the LQR feedback control law for (\ref{eq:comb w dyn}) to get the new dynamics.  The result is a set of equations of the same form as (\ref{eq:comb dyn 1})-(\ref{eq:comb dyn 3}) but with $B^{LQR}_{\left(\cdot\right)}$ instead of $B_{\left(\cdot\right)}$ (and analogous calculations for the new cost function).  The expressions for the various $B$ matrices and cost function components are straightforward to compute (i.e., using limits of the kind expressed in (\ref{eq:y lim 1})-(\ref{eq:y lim 2})) but are omitted due to their length.

For the purposes of optimal control, the time scale separation can be dealt with in one of two ways: a) leave it in the $\epsilon > 0$ case and explicitly solve for all of the $\psi_y$ and $\psi_w$ trajectories, or b) use the $\epsilon \rightarrow 0$ limit to collapse everything into the slow time scale, thereby significantly reducing the size of the problem at the cost of potentially reducing the accuracy of the solution.  Once that choice has been made, we can then choose whether to use the original actuator dynamics or the KO-based LQR dynamics in our optimization.

\subsection{Analytical Tools}

The KO provides opportunities to use a wide range of analytical tools for linear systems, both controlled and uncontrolled.  For the purpose of this paper, however, we will focus on three.  For a dynamical system $x_{t+1} = A x_t$, we can consider the maximum initial growth of a transient, $\log \left\| A\right\|_2$, as well as the maximum transient growth lower bound and complex stability radius, which are, respectively, \cite{trefethen2005spectra}

\begin{gather}
		\sup \limits_{\left|z\right| > 1} \left(\left|z\right| - 1\right) \left\| \left(zI - A\right)^{-1}\right\|_2 \\
		\left[ \sup \limits_{\left|z\right| = 1} \left\| \left(zI - A\right)^{-1}\right\|_2\right]^{-1}
\end{gather}

Long-term stability, measured by matrix eigenvalues, is important, but as Trefethen and Embree discuss, large transient growth can make a stable system just as unusable, in practice, as an unstable system \cite{trefethen2005spectra}.  We will show how the interactions between time scales affects these measures of stability by considering how they change when applied to the slow time scale dynamics with and without the $\epsilon \rightarrow 0$ limit.  For example, for the original time scale separated case, we will consider the difference between evaluating these quantities when $A = K_{xx}$ and when $A = K_{comb,xx}$.  The latter encapsulates the feedback between the two time scales while the former does not.

\onecolumn
\par\noindent\rule{\dimexpr(0.5\textwidth-0.5\columnsep-0.4pt)}{0.4pt}%
\rule{0.4pt}{6pt}

\begin{align}
		\psi_x \left(x_{t+1}\right) &= K_{xx} \left(\psi_x \left(x_t\right) - K_{xw} \frac{1}{m} \sum \limits_{n=0}^{m-1} \psi_w \left(w_{t + n \tau} \right) - K_{xy} \frac{1}{m} \sum \limits_{n=0}^{m-1} \psi_y \left(y_{t + n\tau}\right) \right) \nonumber \\
		&+ K_{xw} \frac{1}{m} \sum \limits_{n=0}^{m-1} \psi_w \left(w_{t + n \tau} \right) + K_{xy} \frac{1}{m} \sum \limits_{n=0}^{m-1} \psi_y \left(y_{t + n\tau}\right) \label{eq:comb x dyn} \\
		\psi_y \left(y_{t+\left(n+1\right)\tau}\right) &= K_{yy} \left(\psi_y \left(y_{t + n \tau}\right) - K_{yx} \psi_x \left(x_t\right) - K_{yw} \psi_w \left(w_{t + n \tau}\right) \right) + K_{yx} \psi_x \left(x_t\right) + K_{yw} \psi_w \left(w_{t + n \tau}\right) \label{eq:comb y dyn} \\
		\psi_w \left(w_{t+\left(n+1\right)\tau}\right) &= K_{ww} \left(\psi_w \left(w_{t + n \tau}\right) - K_{wx} \psi_x \left(x_t\right) - K_{xy} \psi_y \left(y_{t + n \tau}\right) - K_{wu} \psi_u \left(u_t\right) \right) \nonumber \\
		& + K_{wx} \psi_x \left(x_t\right) + K_{xy} \psi_y \left(y_{t + n \tau}\right) + K_{wu} \psi_u \left(u_t\right) \label{eq:comb w dyn} \\
		\psi_{y,t}^T P \psi_{y,t} + p^T \psi_{y,t} &= \min_{\psi_{w,t}} \left[ \sum_{i,j} \psi_{i,t}^T Q_{ij} \psi_{j,t} + \sum_i c_i^T \psi_{i,t} + \psi_{y,t+1}^T P \psi_{y,t+1} + p^T \psi_{y,t+1} \right] \label{eq:comb bellman} \\
		P &= K_{yy}^T P K_{yy} - \frac{1}{2} F^T Q F + Q_{yy} \\
		\left(I - K_{yy}\right)^T p &= - F^T q + \left[ Q_{yx} + Q_{xy}^T + 2 K_{yy}^T P \left(I-K_{yy}\right) K_{yx} \right] \psi_{x,t} + \left(Q_{uy}^T + Q_{yu} \right) \psi_{u,t} + c_y \\
		\frac{1}{2} d^T Q d &= \psi_{x,t}^T Q_{xx} \psi_{x,t} + \psi_{u,t}^T Q_{uu} \psi_{u,t} + \psi_{x,t}^T \left(Q_{xu} + Q_{ux}^T \right) \psi_{u,t} + \left[c_x^T + p^T \left(I - K_{yy}\right) K_{yx}\right] \psi_{x,t} \nonumber \\
		&+ c_u^T \psi_{u,t} + c_0 \label{eq:comb lqr end}
\end{align}

\par\noindent\rule{\dimexpr(0.5\textwidth-0.5\columnsep-0.4pt)}{0.4pt}%
\rule[-6pt]{0.4pt}{6pt}
\begin{multicols}{2}

\section{Computational Demonstration}

\subsection{Problem Formulations}

Let us assume that we have a van der Pol oscillator operating at a slow time scale ($\mathbf{x}$) connected to a Duffing oscillator operating at a fast time scale ($\mathbf{y}$), where $\epsilon = 100$, and that we have an actuator $w$ evolving at the fast time scale according to a PI loop trying to match $y_1$ to the control input $u$:

\begin{align}
		\dot{\mathbf{y}} &= \frac{1}{\epsilon} \left\{ \begin{array}{c} y_2 \\ -2 y_2 - y_1 - y_1^3 + 0.5 x_2^2 - 2w \end{array} \right\} \\
		\dot{\mathbf{x}} &= \left\{ \begin{array}{c} x_2 \\ -0.5 \left(1 - x_1^2\right) x_2 - x_1 + y_1 \end{array} \right\} \\
		w \left(t\right) &= K_1 \left(y_1 - u\right) + K_2 \int \limits_0^t \left(y_1 - u\right) dt \\
		\Rightarrow \dot{w} &=  K_1 \dot{y}_1  + K_1 \dot{u} + K_2 \left(y_1 - u\right) \nonumber \\ 
		&= \frac{1}{\epsilon} K_1 y_2 + K_1 \dot{u} + K_2 \left(y_1 - u\right) 
\end{align}

The coupling functions $0.5 x_2^2$ and $y_1$ mean that the van der Pol kinetic energy affects the Duffing, the Duffing position affects the van der Pol oscillator (acting like a force), and the actuator input acts like a force on the Duffing.  If $K_1 \approx K_2$, then $\dot{w} \approx \frac{1}{\epsilon} K_1 y_2$, which is really just proportional control.  If $K_1 \approx K_2 / \epsilon$, then $K_1 \dot{u} \approx 0$ and we retain PI control.  Based on some unoptimized manual tuning, we chose $K_1 \epsilon =1$ and $K_2=1.0$.  For the system-level optimal control, let us assume that we have a running cost of $y_1^2 + y_2^2 + x_1^2 + x_2^2 + w^2$.  For LQR control (replacing PI control), let us assume a running cost of $\left(y_1 - u\right)^2 + w^2$.  We also constrain $x$ and $u$ to lie between -1 and 1; the system-level objective is unchanged.  When we are only considering the time scale separated case, we remove the actuator equations and variables from this formulation, and when we are only considering the hierarchical control case, we remove the slow-scale equations and variables but keep the control input $u$ as a constant to be optimized (i.e., we solely consider the fast time scale dynamics).  Also, in the full problem, $w$ and $x$ do not appear directly in each other's dynamics, so we can use this domain information to set $K_{xw}$ and $K_{wx}$ to be identically zero (see \cite{bakker2022deception} for more on this kind of domain knowledge incorporation).

\subsection{Koopman Implementation}

To learn a KO model for each case, we used state-inclusive observables for $\psi_x$, $\psi_y$, and $\psi_w$ with the nonlinear observables calculated using small neural networks.  We used 12 nonlinear observables for $\psi_x$ and $\psi_y$ and four nonlinear observables for $\psi_w$; $u$ was kept unlifted, as is common in the KO literature (e.g., \cite{korda2018linear,king2021solving}).  To train the neural networks (implemented using Keras in Tensorflow), we chose 1e4 random initial conditions in the $\left[-1,1\right]^n$ hypercube and simulated the system $\Delta t = 0.1$ seconds forward in time to produce our training data.  For the slow time scale data, we only recorded data at $t = 0$ and $t = \Delta t$, but for the fast time scale data, we recorded data at $\Delta t/\epsilon$ time steps during that interval; in the hierarchical case, with only the fast dynamics present, only the first $\Delta t/\epsilon$ step was used for training.  The training function consisted of mean squared prediction error and stability enforcement terms in all three cases.  Stability enforcement involved calculating the eigenvalues of the relevant matrices and penalizing eigenvalues with magnitudes greater than one.  It was applied both to all of the relevant matrices in the $\epsilon > 0$ case (e.g., $K_{xx}$, $K_{yy}$) and to the slow time scale matrices produced in the $\epsilon \rightarrow 0$ limit (e.g., $K_{comb,xx}$, $B_{xx}$).  The cases with hierarchical control also included terms for the solution of the LQR equations (e.g., (\ref{eq:comb lqr end})) and for the internal consistency of the calculated $\psi_w$ in the LQR feedback law (as discussed in Section \ref{sec:hierarchical control}).  The optimal control problems were formulated in Pyomo \cite{hart2017} and solved using IPOPT \cite{wachter2006implementation}.

\section{Results}

\subsection{Time Scale Separation}

\begin{center}
\begin{table}[H]
    \centering
    \caption{Mean RMS Koopman Prediction Error over 100 Random Trajectories, 100 Time Steps, Time Scale Separation}
    \label{tab:rms tss}
    \vspace{6pt}
    \begin{tabular}{c|cc}
        Dynamics &$\epsilon \neq 0$ Error &$ \epsilon \rightarrow 0$ Error \\ 
        \hline
        Slow ($x$) &0.111 &0.111 \\
        Fast ($y$) &-- &0.036
    \end{tabular}
\end{table}
\end{center}

\begin{center}
\begin{table}[H]
    \centering
    \caption{Stability Results, Time Scale Separation}
    \label{tab:stability tss}
    \vspace{6pt}
    \begin{tabular}{c|cc}
        Stability Measure &$K_{xx}$ &$K_{comb,xx}$ \\ 
        \hline
        Maximum Initial Transient Growth &0.419 &0.484 \\
        Maximum Transient Growth Lower Bound &3.524 &7.844 \\
				Complex Stability Radius &4.25e-3 &8.09e-4
    \end{tabular}
\end{table}
\end{center}

Table \ref{tab:rms tss} indicates fairly good predictive accuracy for our time scale separated model, but that accuracy was noticeably better at the fast scale than the slow one; the slow scale prediction relies on many repeated evaluations of the fast scale predictions, so the fast scale errors compounded.  However, the $\epsilon \rightarrow 0$ limit model (i.e., $K_{comb,xx}$) was just as accurate as the $\epsilon \neq 0$ model, so we gained meaningful information from that derivation.  Table \ref{tab:stability tss}, however, indicates that the feedback between scales has a destabilizing effect on the system.  In particular, the lower bound on the maximum transient growth is significantly greater for $K_{comb,xx}$ than for $K_{xx}$, which means that cross-scale feedbacks have the potential to temporarily amplify transients before the system settles down.  The stability radii are also very small in both cases: with or without the cross scale feedbacks, the system is close to being unstable in the domain over which we trained our model.

\subsection{Hierarchical Control}

We could consider predictive accuracy for the hierarchical control model as we did for the time scale separated model, but since the point of the hierarchical model is to calculate optimal control policies, we instead focus on the optimal control results.  We considered 100 different random initial conditions and solved both the LQR and PI formulations from those starting points.  The LQR optimal control produced a solution with a median relative improvement (measured using the cost function) of 48.2\% compared to the PI optimal control, and the difference in computational costs were negligible (each taking roughly 0.1 seconds to solve).  Because the key optimization parameter $u$ was just a single scalar, we were able to exhaustively sample over its range, and we found that the LQR optimal solution was almost identical to the true optimal LQR solution (a median relative increase of 0.001\%) and the PI optimal solution was also very close to true optimal PI solution (a median relative increase of 0.14\%); differences are essentially due to predictive errors in the approximate Koopman model.

\subsection{Hierarchical Control with Time Scale Separation}

\begin{center}
\begin{table}[H]
    \centering
    \caption{Mean RMS Koopman Prediction Error over 100 Random Trajectories with Constant Control Inputs, 100 Steps, Time Scale Separation with Hierarchical Control}
    \label{tab:rms comb}
    \vspace{6pt}
    \begin{tabular}{c|cc}
        Dynamics &$\epsilon \neq 0$ Error &$ \epsilon \rightarrow 0$ Error \\ 
        \hline
        Slow ($x$) &0.162 &0.115 \\
        Fast ($y$) &-- &0.053 \\
				Fast ($w$) &-- &0.092
    \end{tabular}
\end{table}
\end{center}

Table \ref{tab:rms comb} indicates predictive accuracy comparable to that of Table \ref{tab:rms tss}; adding controls made the system more difficult to model.  Again, predictive accuracy was better at the faster time scales, but interestingly, the $\epsilon \rightarrow 0$ predictive error was better than the $\epsilon \neq 0$ predictive error: applying the $\epsilon \rightarrow 0$ limit actually allowed the model to avoid some of the compounded error in the fast time scale predictions.

\begin{center}
\begin{table}[H]
    \centering
    \caption{Stability Results, Time Scale Separation with Hierarchical Control}
    \label{tab:stability comb}
    \vspace{6pt}
    \begin{tabular}{c|ccc}
        Stability Measure &$K_{xx}$ &$B_{xx}$ &$B^{LQR}_{xx}$ \\ 
        \hline
        Maximum Initial Transient Growth &0.767 &0.79 &0.81\\
        Maximum Transient Growth Bound &7.18 &7.03 &6.89\\
				Complex Stability Radius &2.23e-3 &1.31e-3 &1.33e-3
    \end{tabular}
\end{table}
\end{center}

Table \ref{tab:stability comb} presents  the system stability properties and indicates how those stability properties could be affected by introducing LQR control.  In this case, the differences were not meaningful -- both between $\epsilon \rightarrow 0$ and $\epsilon \neq 0$ as well as between PI and LQR control at the fast time scale.  We might expect LQR control to increase the overall stability of the system by converging more quickly to the fast time scale equilibrium, but that did not necessarily translate into greater stability at the slow time scale.

Again, we considered 100 different random initial conditions and solved both the LQR and PI formulations from those starting points.  In order to make the problem computationally feasible, we used the $\epsilon \rightarrow 0$ approximation.  This reduced the problem size by orders of magnitude (both the number of variables and the number of constraints), and since computational cost scales superlinearly with problem size, the computational cost savings would have been even more significant.  As it was, the optimization solved in fractions of a second (for both LQR and PI).  Furthermore, the greater predictive accuracy of the $\epsilon \rightarrow 0$ model suggests that it might work better for the optimal control calculations in any case; in evaluating the calculated control policies, we simply initialized the original system dynamics with random $w$ and $y$ initial conditions.  On average, the LQR optimal control produced solutions that were very similar, in terms of overall cost, to the PI optimal control (typically within a few percent, and not always better).  Unlike the hierarchical control case, we cannot exhaustively sample over the possible control policies, but if we assume a constant control policy (e.g., as might be produced in a rule-based supervisory control system), we can find an optimal constant policy against which to compare our solutions.  In doing this, we found that the median improvements for the LQR and PI optimal solutions, compared to the constant solutions, were almost 40\% in both cases.  It did not matter much whether PI or LQR control was being used, but optimizing a time-varying supervisory control policy did matter.

\section{Discussion and Conclusions}

In this paper, we used the KO to model systems that have time scale separated dynamics, hierarchical control, and both.  We then showed how those KO models could be used to calculate optimal control policies (at multiple levels) and to gain insight into the impacts of cross-scale feedbacks on stability; this paper simply provided a few examples of how the KO could be used for such systems.  The particular structure of the KO models presented here facilitated the analysis and control calculations.  Separating the actuators and the fast time scale states allowed us to replace the PI local control with LQR control without having the relearn the model, for example; the LQR policy did have to be learned during the model training process, however.  Similarly, separating the KO representations of the fast and slow time scales enabled the calculation of the $\epsilon \rightarrow 0$ approximate dynamics (both with and without control), and that in turn supported computationally efficient optimal control solutions based on the $\epsilon \rightarrow 0$ dynamics.

Substituting LQR for PI control was very beneficial in the hierarchical control case but had little effect when combined with time scale separated dynamics.  In hindsight, this makes sense.  The LQR and PI control laws focus on getting the fast time scale dynamics to converge to their equilibrium point.  Once there, the system behaves essentially like a DAE -- changes to $y$ become dominated by $x$'s dynamics (and how those dynamics slowly change $y$'s equilibrium point), so using LQR vs. PI does not make much difference.  This suggests an alternative to either PI or LQR control: design a feedback control policy with the slow time scale in mind.  For example, we could choose $F$ and $d$ (see Section \ref{sec:hierarchical control}) such that they maximize the stability of the $\epsilon \rightarrow 0$ dynamics (i.e., minimize the largest eigenvalue of $B_{xx}$) or minimize transient growth (i.e., $\log \left\| B_{xx}\right\|_2$) while still keeping the fast time scale dynamics stable.  Alternatively, we could focus on the controllability of the $\epsilon \rightarrow 0$ system as measured by the controllability gramian calculated from $B_{xx}$ and $B_{xu}$ (e.g., see \cite{sinha2022data}).  The construction and use of such control policies would be relatively straightforward and would also constitute a potentially fruitful avenue for future research.  Other areas for future research include applications to larger and more realistic systems or to systems with more than two different time scales of operation.



\bibliographystyle{ieeetr}
\bibliography{bib}

\end{multicols}

\end{document}